\documentclass[submission%
]{dmtcs}

\usepackage[latin1]{inputenc}
\usepackage{subfigure}

\usepackage{amsmath}

\newtheorem{defn}{Definition}
\newtheorem{thm}{Theorem}

\newtheorem{lemma}{Lemma}
\newtheorem{cor}{Corollary}

\usepackage[font=small, labelfont=bf, width=0.8\textwidth]{caption}

\DeclareMathOperator{\diag}{diag}

\DeclareMathOperator{\Tr}{Tr}

\def\be{\begin{equation}}
\def\ee{\end{equation}}

\def\phid{\phi^\dagger}

\author{Luigi Cantini\addressmark{1}\thanks{Email: \email{luigi.cantini@u-cergy.fr}}, Jan de Gier\addressmark{2}\thanks{Presenting author, Email: \email{jdgier@unimelb.edu.au}}
\and Michael Wheeler\addressmark{2}\thanks{Email: \email{wheelerm@unimelb.edu.au}}}
\title{Matrix product and sum rule for Macdonald polynomials}
\address{\addressmark{1}Laboratoire de Physique Th\'eorique et Mod\'elisation (CNRS UMR 8089), Universit\'e de Cergy-Pontoise, F-95302 Cergy-Pontoise, France\\ \addressmark{2}ARC Centre of Excellence for Mathematical and Statistical Frontiers (ACEMS), School of Mathematics and Statistics, The University of Melbourne, VIC 3010, Australia}
\keywords{Macdonald polynomials, matrix product ansatz, Yang--Baxter equation, asymmetric exclusion process}
\begin{document}
\maketitle
\begin{abstract}
\paragraph{Abstract.} 
We present a new, explicit sum formula for symmetric Macdonald polynomials $P_\lambda$ and show that they can be written as a trace over a product of (infinite dimensional) matrices. These matrices satisfy the Zamolodchikov--Faddeev (ZF) algebra. We construct solutions of the ZF algebra from a rank-reduced version of the Yang--Baxter algebra. As a corollary, we find that the normalization of the stationary measure of the multi-species asymmetric exclusion process is a Macdonald polynomial with all variables set equal to one.

\paragraph{R\'esum\'e.}
Nous pr\'esentons une nouvelle formule explicite pour les polyn\^{o}mes sym\'etriques de Macdonald $P_{\lambda}$, et d\'emontrons qu'ils peuvent \^{e}tre \'ecrits comme une trace d'un produit de matrices de dimension infinie. Ces matrices forment une repr\'esentation de l'alg\`ebre de Zamolodchikov--Faddeev (ZF). Nous \'elaborons des solutions de l'alg\`ebre ZF en utilisant une version de l'alg\`ebre Yang--Baxter de rang r\'eduit. Comme corollaire, nous constatons que la normalisation de la mesure stationnaire du processus d'exclusion asym\'etrique est un polyn\^{o}me de Macdonald avec toutes les variables rendues \'egale \`a un.     

\paragraph{Samenvatting.} 
We presenteren een nieuwe, expliciete som formule voor symmetrische Macdonald polynomen $P_\lambda$ en laten zien dat deze geschreven kunnen worden als een spoor over een produkt van (oneindig dimensionale) matrices. Deze matrices voldoen aan de Zamolodchikov--Faddeev (ZF) algebra. We construeren oplossingen van de ZF algebra uit een rank-gereduceerde versie van de Yang-baxter algebra. Een gevolgtrekking is dat de normalisatie van de stationaire toestand van het asymmetrische exclusie process met meerdere soorten deeltjes gegeven wordt als een Macdonald polynoom waarin alle variabelen gelijk aan \'e\'en worden gezet.
\end{abstract}

\section{Introduction}
\label{sec:in}
Symmetric Macdonald polynomials \cite{macd88,MacdBook} are a family of multivariable orthogonal polynomials indexed by partitions, whose coefficients depend rationally on two parameters $q$ and $t$. In the case $q=t$ they degenerate to the celebrated Schur polynomials, which are central in the representation theory of both the general linear and symmetric groups. Let $m_\lambda$ denote the monomial symmetric polynomial indexed by a partition $\lambda$, i.e. the symmetric polynomial defined as the sum of all monomials $x^\mu=x_1^{\mu_1}\cdots x_n^{\mu_n}$ where $\mu$ ranges over all distinct permutations of $\lambda=(\lambda_1,\ldots,\lambda_n)$. The Macdonald polynomials are defined as follows:
\begin{defn}Let $\langle \cdot, \cdot \rangle$ denote the Macdonald inner product on power sum symmetric functions (\cite{MacdBook}, Chapter VI, Equation (1.5)),
where $<$ denotes the dominance order on partitions (\cite{MacdBook}, Chapter I, Section 1). The Macdonald polynomial $P_{\lambda}(x_1,\dots,x_n;q,t)$ is the unique homogeneous symmetric polynomial in $(x_1,\dots,x_n)$ which satisfies
\begin{align*}
\langle P_{\lambda}, P_{\mu} \rangle &= 0, \ \lambda \neq \mu,\\
P_{\lambda}(x_1,\dots,x_n;q,t)
&=
m_{\lambda}(x_1,\dots,x_n)
+
\sum_{\mu < \lambda}
c_{\lambda,\mu}(q,t)
m_{\mu}(x_1,\dots,x_n),
\end{align*}
i.e. the coefficients $c_{\lambda,\mu}(q,t)$ of the lower degree terms are completely determined by the orthogonality conditions.
\end{defn}

Up to normalization, Macdonald polynomials can alternatively be defined as the unique eigenfunctions of certain linear difference operators acting on the space of all symmetric polynomials \cite{MacdBook}. They can also be expressed combinatorially as multivariable generating functions \cite{HaglundHL1,HaglundHL2,RamY}, or via symmetrization of non-symmetric Macdonald polynomials that are computed from Yang--Baxter graphs \cite{Lasc01,Lasc}.

The purpose of this article is to report on an explicit matrix product formula for Macdonald polynomials \cite{CantinidGW} inspired by recent results on the multi-species asymmetric exclusion process, and to provide an explicit sum rule for calculating Macdonald polynomials resulting from the matrix product formula \cite{dGW}.

In the following we need the polynomial representations of the Hecke algebra of type $A_{n-1}$, with generators $T_{i}$ given by 
\begin{align}
\label{hecke-gen}
T_i
=
t-\frac{tx_i-x_{i+1}}{x_i-x_{i+1}}(1-s_i),
\qquad
1 \leq i \leq n-1,
\end{align} 
where $s_i$ is the transposition operator with action $s_i f(\dots,x_i,x_{i+1},\dots) = f(\dots,x_{i+1},x_i,\dots)$ on functions in $(x_1,\dots,x_n)$. It can be verified that the operators \eqref{hecke-gen} indeed give a faithful representation of the Hecke algebra:
\begin{align*}
(T_{i}-t)(T_{i}+1)=0,
\qquad
T_{i} T_{i\pm 1} T_{i} = T_{i \pm 1} T_i T_{i\pm 1},
\qquad
T_i T_j = T_j T_i, \ |i-j| > 1.
\end{align*}
In view of the relations for the generators, we can define $T_{\sigma}$ unambiguously as any product of simple transpositions $T_{i}$ which gives the permutation $\sigma$. 

\section{Main results}
\label{sec:main}
\begin{defn}[mASEP basis]
\label{def:mASEP}
For all partitions $\lambda$, let $\{f_{\mu=\sigma\circ\lambda}\}_{\sigma\in S_n}$ be a set of homogeneous degree $|\lambda|$ polynomials in $n$ variables, defined by
\begin{align*}
&T_i f_{\mu_1,\dots,\mu_n}(x_1,\dots,x_n;q,t)
=
f_{\mu_1,\dots,\mu_{i+1},\mu_i,\dots,\mu_n}
(x_1,\dots,x_n;q,t),
\qquad
\text{when}\ \ 
\mu_i > \mu_{i+1},\\
&T_i f_{\mu_1,\dots,\mu_n}(x_1,\dots,x_n;q,t)
=
t f_{\mu_1,\dots,\mu_{i+1},\mu_i,\dots,\mu_n}
(x_1,\dots,x_n;q,t),
\qquad
\text{when}\ \ 
\mu_i = \mu_{i+1},\\
&f_{\mu_n,\mu_1,\dots,\mu_{n-1}}(qx_n,x_1,\ldots,x_{n-1};q,t) = q^{\mu_n}f_{\mu_1,\dots,\mu_n}(x_1,\dots,x_n;q,t).
\end{align*}
\end{defn}
As explained in \cite{CantinidGW}, the polynomials $\{f_{\mu}\}$ are related to the non-symmetric Macdonald polynomials \cite{Chera,Cherb,Opdam} via an invertible triangular change of basis, and hence form a basis for the ring of polynomials in $n$ variables. Specializations of these polynomials at $q=x_1=\cdots=x_n=1$ give stationary particle configuration probabilities of the multi-species asymmetric exclusion process on a ring. 

Summing over all $\{f_{\mu=\sigma\circ\lambda}\}_{\sigma\in S_n}$ results in a symmetric Macdonald polynomial \cite{macd95}:

\begin{lemma}
\label{le:symmetry}
Let $\lambda$ be a partition. Then
\[
P_{\lambda}(x_1,\ldots,x_n;q,t) = \sum_{\sigma \in S_n} f_{\sigma \circ \lambda}(x_1,\dots,x_n;q,t).
\] 
\end{lemma}

Suppose we have (semi-infinite) matrices $A_0(x), A_1(x),\ldots, A_r(x)$ and $S$ satisfying the following exchange relations
\begin{align}
A_i(x)A_i(y) 
&= 
A_i(y)A_i(x),
\label{Exchange0}
\\
t A_j(x)A_i(y)
-
\frac{tx-y}{x-y}
\Big( A_j(x)A_i(y) - A_j(y)A_i(x) \Big)
&= A_i(x) A_j(y),
\label{Exchange1}
\\
S A_i(q x) 
&= 
q^{i} A_i(x)S, 
\label{Exchange2}
\end{align}
for all $0 \leq i<j \leq r$. The main result of \cite{CantinidGW} is a matrix product formula for the polynomials $f_{\sigma\circ\lambda}$:
\begin{thm}
\label{th:MP}
There is an explicit representation of $A_0(x), A_1(x),\ldots, A_{r}(x)$ and $S$ with $r=\lambda_1$ satisfying \eqref{Exchange0}, \eqref{Exchange1} and \eqref{Exchange2} such that $f_\mu$ can be written as a matrix product, i.e.
\[
\Omega_{\lambda} f_{\mu}(x_1,\ldots,x_n) = \Tr \Big[ A_{\mu_1}(x_1) \cdots A_{\mu_n}(x_n) S\Big],
\label{MPA}
\]
where $\mu$ is a permutation of $\lambda$ and $\Omega_{\lambda}$ is a normalization factor which only depends on the partition $\lambda$.
\end{thm}

\begin{cor}
It follows from Lemma~\ref{le:symmetry} that the symmetric Macdonald polynomial $P_\lambda$ can be expressed as a sum over matrix product formulas. The specialization of $P_\lambda$ at $q=1$ and $x_i=1$ $(i=1,\ldots,n)$ is the normalization of the stationary state of a multi-species asymmetric exclusion process on a ring.
\end{cor}

As a consequence of Theorem~\ref{th:MP} we derive an explicit sum formula for Macdonald polynomials \cite{dGW}. To formulate this result we need to prepare some notation. Let $\lambda$ be a partition whose largest part is $\lambda_1=r$. For all $0 \leq k \leq r$, we define a partition $\lambda[k]$ by replacing all parts in $\lambda$ of size $\leq k$ with 0. For example, for $\lambda = (3,3,2,1,1,0)$ we have
\begin{align*}
\lambda[0] = (3,3,2,1,1,0),
\quad
\lambda[1] = (3,3,2,0,0,0),
\quad
\lambda[2] = (3,3,0,0,0,0),
\quad
\lambda[3] = (0,0,0,0,0,0),
\end{align*}
and in general $\lambda[0] = \lambda$, $\lambda[r] = 0$.

\begin{thm}[\cite{dGW}]
\label{th:sum}
The Macdonald polynomial $P_{\lambda}$ can be written in the form
\begin{align}
\label{eq:sum}
P_{\lambda}(x_1,\dots,x_n;q,t)
=
\sum_{\sigma \in S_{\lambda} }
T_{\sigma}
\circ
x_{\lambda}
\circ
\prod_{i=1}^{r-1}
\left(
\sum_{\sigma \in S_{\lambda[i]}}
C_i
\left(
\begin{array}{@{}c@{}}
\lambda[i-1]
\medskip \\ 
\sigma\circ\lambda[i]
\end{array}
\right)
T_{\sigma}
\circ
x_{\lambda[i]}
\circ
\right)
1
\end{align}
with coefficients\footnote{We will use three notations for the coefficients interchangeably:
$
C_i( \lambda, \mu) 
\equiv
C_i\left(\begin{array}{@{}c@{}} \lambda \\ \mu \end{array}\right)
\equiv
C_i\left(
\begin{array}{@{}c@{}c@{}c@{}} \lambda_1 & \cdots & \lambda_n 
\\ 
\mu_1 & \cdots & \mu_n \end{array}
\right)
$.} that satisfy $C_i( \lambda, \mu) = 0$ if any $0 < \lambda_k < \mu_k$, and otherwise
\begin{align*}
\label{coeff}
&
C_i(\lambda,\mu)
\equiv
C_i\left(\begin{array}{@{}c@{}c@{}c@{}} 
\lambda_1 & \cdots & \lambda_n \\ 
\mu_1 & \cdots & \mu_n \end{array}\right)
=
\prod_{j = i+1}^{r}
\left(
q^{(j-i)a_j(\lambda,\mu)}
\prod_{k=1}^{b_j(\lambda,\mu)}
\frac{1-t^k}{1-q^{j-i}t^{\lambda'_i-\lambda'_j+k}}
\right),
\end{align*}
where
$\ a_i(\lambda,\mu)
=
\#\{ (\lambda_k,\mu_k): \lambda_k = 0,\ \mu_k = i \}\ $ and
$\ 
b_i(\lambda,\mu)
=
\#\{ (\lambda_k,\mu_k): i = \lambda_k > \mu_k \}.
$
\end{thm}

We point out that the formula given in \eqref{eq:sum} has many structural features in common with the work of Kirillov and Noumi \cite{Kirillov_Noumi1,Kirillov_Noumi2}. In these papers the authors construct families of raising operators, which act on Macdonald polynomials by adding columns to the indexing Young diagram. In \cite{Kirillov_Noumi1} the raising operators have an analogous form to Macdonald $q$-difference operators, while in \cite{Kirillov_Noumi2} the raising operators are constructed in terms of generators of the affine Hecke algebra. In both papers the Macdonald polynomial is obtained by the successive action of such raising operators on $1$, the initial state. It would be very interesting to find a precise connection between the results of \cite{Kirillov_Noumi1,Kirillov_Noumi2} and our formula \eqref{eq:sum}, if one exists. 

\paragraph{Two variable example.}
Let us demonstrate \eqref{eq:sum} in the case $\lambda = (3,1)$. We have 
\begin{align*}
\lambda[0] = (3,1),\ \lambda[1] = (3,0),\ \lambda[2] = (3,0),\ \lambda[3] = (0,0),
\end{align*}
and $S_{\lambda[0]} = S_{\lambda[1]} = S_{\lambda[2]} = S_2$. Hence
\begin{multline*}
P_{(3,1)}(x_1,x_2;q,t)
=
\sum_{\sigma \in S_2}
T_{\sigma}
\circ
x_1 x_2
\circ
\sum_{\rho \in S_2}
C_1\left( 
\begin{array}{@{}l@{}l@{}}
\lambda[0]_1 & \lambda[0]_2 \\
\lambda[1]_{\rho_1} & \lambda[1]_{\rho_2} 
\end{array}\right) \times
\\
T_{\rho}
\circ
x_1
\circ
\sum_{\pi \in S_2}
C_2\left( 
\begin{array}{@{}l@{}l@{}}
\lambda[1]_1 & \lambda[1]_2 \\
\lambda[2]_{\pi_1} & \lambda[2]_{\pi_2} 
\end{array}\right)
T_{\pi}
\circ
x_1.
\end{multline*}
We compute each sum in turn, starting with the rightmost:
\begin{align*}
\sum_{\pi \in S_2}
C_2\left( 
\begin{array}{@{}l@{}l@{}}
\lambda[1]_1 & \lambda[1]_2 \\
\lambda[2]_{\pi_1} & \lambda[2]_{\pi_2} 
\end{array}\right)
T_{\pi}
\circ
x_1
&=
x_1
+
C_2\left( 
\begin{array}{@{}cc@{}}
3 & 0 \\
0 & 3
\end{array}\right)
T_{1}
\circ
x_1
=
x_1
+
q \left( \frac{1-t}{1-qt} \right)
x_2.
\end{align*}
Combining with the middle sum, we find that
\begin{align*}
\sum_{\rho \in S_2}
C_1\left( 
\begin{array}{@{}l@{}l@{}}
\lambda[0]_1 & \lambda[0]_2 \\
\lambda[1]_{\rho_1} & \lambda[1]_{\rho_2} 
\end{array}\right)
T_{\rho}
\circ
\left(
x_1^2
+
q\frac{1-t}{1-qt} x_1 x_2
\right)
&=
x_1^2
+
q\frac{1-t}{1-qt} x_1 x_2,
\end{align*}
where we have used the fact that
$C_1\left( 
\begin{array}{@{}cc@{}}
3 & 1 \\
3 & 0
\end{array}\right)= 1$ and 
$C_1\left( 
\begin{array}{@{}cc@{}}
3 & 1 \\
0 & 3
\end{array}\right)= 0$. Combining everything and passing to the leftmost sum, we have
\begin{align*}
P_{(3,1)}(x_1,x_2;q,t)
&=
\sum_{\sigma \in S_2}
T_{\sigma}
\circ
\left(
x_1^3 x_2
+
q  \frac{1-t}{1-qt} 
x_1^2 x_2^2
\right)
\\
&=
(1+T_{1})\circ
\left(
x_1^3 x_2
+
q  \frac{1-t}{1-qt} 
x_1^2 x_2^2
\right)
=
x_1^3 x_2
+
\frac{1-t+q-qt}{1-qt}
x_1^2 x_2^2
+
x_1 x_2^3.
\end{align*}

\section{Yang--Baxter and Zamolodchikov--Faddeev algebras}

In the remainder of this note we sketch the proof of Theorem~\ref{th:MP}. It is relatively simple to show that if the functions $f_\mu$ defined in Definition~\ref{def:mASEP}, with $\mu$ a permutation of a partition $\lambda$, can be written as a matrix product formula as in Theorem~\ref{th:MP}, then the matrices need to satisfy the exchange relations \eqref{Exchange0}--\eqref{Exchange2}. These can be conveniently rewritten.

As before, let $r=\lambda_1$ be the largest part of $\lambda$, we call $r$ the rank. It was shown in \cite{SasaW97} for $r=1$ and general $r$ in \cite{CrampeRV}, see also \cite{CantinidGW}, that \eqref{Exchange0}--\eqref{Exchange2} are equivalent to the Zamolodchikov--Faddeev (ZF) algebra \cite{ZZ1979,Fad1980},
\be
\check{R}(x,y )\cdot\left [\mathbb{A}(x)\otimes \mathbb{A}(y)\right] 
= 
\left [\mathbb{A}(y)\otimes \mathbb{A}(x)\right] ,
\label{eq:ZFdef}
\ee
where $\check{R}$ is a twisted version of the $U_{t^{1/2}}(sl_{r+1})$ $R$-matrix 
and  $\mathbb{A}=\mathbb{A}^{(r)}(x)$ is an $(r+1)$-dimensional operator valued column vector given by
\[
\mathbb{A}^{(r)}(x) = (A_0(x),\ldots, A_r(x))^T.
\]
Let furthermore $E^{(ij)}$ denote the elementary $(r+1) \times (r+1)$ matrix with a single non-zero entry 1 at position $(i,j)$. Then equation \eqref{Exchange2} can be rewritten as
\[
S \mathbb{A}(qx) =q^{\sum_i iE^{(ii)}}\mathbb{A}(x) S,
\]
where the rank $r$ is again implicit, i.e. $\mathbb{A}=\mathbb{A}^{(r)}(x)$ and $S=S^{(r)}$, and the $R$-matrix is explicitly given by
\begin{align*}
\check{R}^{(r)} (x,y)
=
\sum_{i=1}^{r+1}
E^{(ii)}
\otimes
E^{(ii)}
+
\frac{x - y}{t x - y}
\sum_{1 \leq i < j \leq r+1}
\Big(
t E^{(ij)}
\otimes
E^{(ji)}
+
E^{(ji)}
\otimes
E^{(ij)}
\Big)\nonumber
\\
\frac{t-1}{t x - y}
\sum_{1 \leq i < j \leq r+1}
\Big(
x
E^{(ii)}
\otimes
E^{(jj)}
+
y
E^{(jj)}
\otimes
E^{(ii)}
\Big).
\end{align*}

\paragraph{Example of a rank 1 solution to ZF algebra.} 
We give a simple explicit example for the case $r=1$. Using the following functions
\begin{align*}
b^+  &= \displaystyle \frac{t(x-y)}{tx-y},  &  b^-  &= t^{-1} b^+ = \frac{x-y}{tx-y}, \nonumber \\
\\[-\baselineskip]
c^+ & =1-b^+ = \frac{y (t-1)}{tx-y}, & c^- &=1-b^- = \frac{x (t-1)}{tx-y}. \nonumber
\end{align*}
equation \eqref{eq:ZFdef} for $\mathbb{A}(x)=\begin{pmatrix}  1 \\ x \end{pmatrix}$ and $r=1$ explicitly becomes: 
\begin{align}
\label{eq:ZFrank1}
\left(
\begin{array}{cc|cc}
1 & 0 & 0 & 0
\\
0 & c^- &  b^+ & 0
\\
\hline
0 & b^- &c^+ & 0
\\
0 & 0 & 0 & 1
\end{array}
\right)
\cdot
\left[
\left(
\begin{array}{c}
1 \\ x
\end{array}
\right)
\otimes
\left(
\begin{array}{c}
1 \\ y
\end{array}
\right)\right]
=
\left[
\left(
\begin{array}{c}
1 \\ y
\end{array}
\right)
\otimes
\left(
\begin{array}{c}
1 \\ x
\end{array}
\right)
\right].
\end{align}

\subsection{General rank}
For general rank $r$, solutions to \eqref{eq:ZFdef} are more difficult, but can be recovered from the Yang--Baxter algebra  which is given by
\begin{align}
\check{R}(x,y)\cdot\left [L(x)\otimes L(y)\right] 
= 
\left [L(y)\otimes L(x)\right] \cdot \check{R}(x,y),
\label{eq:YBAdef}
\end{align}
where $L(x)=L^{(r)}(x)$ is an $(r+1) \times (r+1)$ operator-valued matrix. The algebra \eqref{eq:YBAdef} is well-studied and many solutions for $L(x)$ are known. For the application to Macdonald polynomials the elements of $L(x)$ are given in terms of generators $\{k,\phi, \phi^{\dagger}\}$ of the $t$-boson algebra:
\begin{align}
\label{t-bos}
\phi k = t k \phi,
\qquad
t \phid k = k \phid,
\qquad
\phi \phid - t \phid \phi = 1-t.
\end{align}

We can construct solutions of \eqref{eq:ZFdef} by rank-reducing the Yang--Baxter algebra \eqref{eq:YBAdef} in the following way. Assume a solution of the following modified Yang--Baxter algebra
\begin{align}
\check{R}^{(r)}(x,y)\cdot\left [\tilde{L}(x)\otimes \tilde{L}(y)\right] 
= 
\left [\tilde{L}(y)\otimes \tilde{L}(x)\right] \cdot \check{R}^{(r-1)}(x,y),
\label{eq:YBAdef_lowrank}
\end{align}
in terms of an $(r+1) \times r$ operator-valued matrix $\tilde{L}(x)=\tilde{L}^{(r)}(x)$ and where the rank of the $R$-matrix on the right hand side is one lower than that on the left hand side. Assume also an operator $s=s^{(r)}$ that satisfies
\[
s\tilde{L}(qx) =q^{\sum iE^{(ii)}}\tilde{L}(x) sq^{-\sum iE^{(ii)}},
\]
Then
\begin{align*}
\mathbb{A}^{(r)}(x) &= \tilde{L}^{(r)}(x)\cdot \tilde{L}^{(r-1)}(x)  \cdots \tilde{L}^{(1)}(x),
\qquad
S^{(r)} = s^{(r)}\cdot s^{(r-1)} \cdots s^{(1)}
\label{eq:nestedTwist}
\end{align*}
gives a solution to \eqref{eq:ZFdef} provided that the operator entries of $\tilde{L}^{(a)}(x)$ commute with those of $\tilde{L}^{(b)}(y)$, for all $a \not=b$. The usual way to ensure this commutativity is to demand that the entries of $\tilde{L}^{(a)}$ act on some vector space $V_a$ while $\tilde{L}^{(b)}$ act on a different vector spaces $V_b$, and indeed we shall adopt this approach. 

\paragraph{Rank 1 example continued.}
The solution to the Yang--Baxter algebra  \eqref{eq:YBAdef} corresponding  to rank $r=1$ is equal to
\[
L^{(1)}(x)=\begin{pmatrix} 1 & \phi \\ x \phi^\dag & x \end{pmatrix},
\]
where the operators $\phi$, $\phi^\dagger$ and $k$ satisfy the $t$-boson relations \eqref{t-bos}. We note that trivialising the $t$-boson by sending $\phi^\dagger,\phi\mapsto 1$ and $k\mapsto 0$, we reduce the rank, and thus obtain the  solution $\mathbb{A}^{(1)}(x)=\tilde{L}^{(1)}(x)$ given in \eqref{eq:ZFrank1}:
\[
\begin{pmatrix} 1 & \phi \\  x \phi^\dag & x \end{pmatrix} \mapsto \begin{pmatrix} 1 & 1 \\ x & x \end{pmatrix}.
\]

\paragraph{Rank 2 solution to ZF algebra.}

The rank $2$ case gives rise to operator valued solutions for $\mathbb{A}^{(2)}(x)$. The associated rank 2 solution to the Yang--Baxter algebra is
\[
L^{(2)}(x)=
\begin{pmatrix}
1 & \phi_1 & \phi_2 \\
x \phi_1^\dag k_2 & x k_2 & 0 \label{eq:osc}\\
x \phi_2^\dag & x \phi_1 \phi_2^\dag & x
\end{pmatrix},
\]
where $\{\phi_1,\phi^{\dag}_1,k_1\}$ and  $\{\phi_2,\phi^{\dag}_2,k_2\}$ are two commuting copies of the $t$-boson algebra \eqref{t-bos}. The map $\phi_1^\dagger,\phi_1\mapsto 1$ and $k_1\mapsto 0$ reduces the rank of $L^{(2)}(x)$ by one
\[
L^{(2)}(x) \mapsto
\begin{pmatrix}
1 & 1 & \phi_2 \\
x k_2  & x k_2 & 0\\
x \phi_2^\dag & x \phi_2^\dag & x
\end{pmatrix}\  \Rightarrow\ 
\tilde{L}^{(2)}(x)=
\begin{pmatrix}
1 &  \phi_2 \\
x k_2 & 0\\
x \phi_2^\dag & x
\end{pmatrix},
\]
where the indices of $t$-bosons are redundant in the final matrix, since we no longer need to distinguish between the two copies of the algebra.

Indeed, we find that \eqref{eq:YBAdef_lowrank} is satisfied
\begin{align*}
\left(
\begin{array}{ccc|ccc|ccc}
1 & 0 & 0 & 0 & 0 & 0 & 0 & 0 & 0
\\
0 & c^- & 0 & b^+ & 0 & 0 & 0 & 0 & 0
\\
0 & 0 &c^-& 0 & 0 & 0 & b^+ & 0 & 0
\\
\hline
0 & b^- & 0 & c^+& 0 & 0 & 0 & 0 & 0
\\
0 & 0 & 0 & 0 & 1 & 0 & 0 & 0 & 0
\\
0 & 0 & 0 & 0 & 0 & c^- & 0 & b^+& 0
\\
\hline
0 & 0 & b^- & 0 & 0 & 0 &  c^+ & 0 & 0
\\
0 & 0 & 0 & 0 & 0 & b^- & 0 & c^+& 0
\\
0 & 0 & 0 & 0 & 0 & 0 & 0 & 0 & 1
\end{array}
\right)\cdot
\left[\left(
\begin{array}{cc}
1 & \phi \\ x k & 0 \\ x \phi^\dag & x
\end{array}
\right)\otimes
\left(
\begin{array}{cc}
1 & \phi \\ y k& 0 \\ y \phi^\dag & y
\end{array}
\right)\right]
= \nonumber
\\
\left[\left(
\begin{array}{cc}
1 & \phi \\ y k  & 0\\ y \phi^\dag & y
\end{array}
\right)\otimes
\left(
\begin{array}{cc}
1 & \phi \\ x k & 0 \\ x \phi^\dag & x
\end{array}
\right)\right]\cdot
\left(
\begin{array}{cc|cc}
1 & 0 & 0 & 0
\\
0 & c^- & b^+& 0
\\
\hline
0 & b^- & c^+ & 0
\\
0 & 0 & 0 & 1
\end{array}
\right).
\end{align*}
We thus construct a solution of the ZF algebra in the following way:
\be
\label{eq:Arank2}
\mathbb{A}^{(2)}(x) = \begin{pmatrix} A_0(x) \\ A_1(x) \\ A_2(x) \end{pmatrix} = \tilde{L}^{(2)}(x)\cdot \tilde{L}^{(1)}(x) = \begin{pmatrix} 1 & \phi \\ x k & 0 \\ x \phi^\dag & x \end{pmatrix}  \begin{pmatrix} 1  \\ x \end{pmatrix} =  \begin{pmatrix} 1  + x \phi \\ k x \\ x \phi^\dag+x^2 \end{pmatrix} ,
\ee
which for $x=1$ is the matrix product solution to the stationary state of the two-species ASEP.

\paragraph{General solution.} For completeness we include the general solution which does not look illuminating in algebraic form. A natural and intuitive combinatorial description of this solution is given in \cite{CantinidGW} but due to lack of space we are not able to present this here. Assume $\lambda\subseteq r^n$ and introduce the following family of $(r-s+2) \times (r-s+1)$ operator-valued matrices $\widetilde{L}^{(s)}(x)$, $1 \leq s \leq r$. We index rows by $i \in \{0,s,\dots,r\}$ and columns by $j \in \{0,s+1,\dots,r\}$\footnote{This unusual indexing of the entries of the matrix is the most convenient for our purposes, since we ultimately want to identify these matrix elements with the partitions $\lambda[s]$ introduced in Section \ref{sec:main}.}, and take
\begin{align*}
\nonumber
\widetilde{L}^{(s)}_{00}
=
1,
\qquad
&
\widetilde{L}^{(s)}_{0j}
=
\phi_j,\ \ 
\text{for }\ s+1 \leq j \leq r,
\qquad
\widetilde{L}^{(s)}_{i0}(x)
=
x \times
\left\{
\begin{array}{ll}
\prod_{l=i+1}^{r}
k_l,
&
i=s
\\
\\
\phi^{\dagger}_i
\prod_{l=i+1}^{r}
k_l,
&
s < i \leq r
\end{array}
\right.
\\
\nonumber
\\
&
\qquad
\widetilde{L}^{(s)}_{ij}(x)
=
x \times
\left\{
\begin{array}{ll}
\prod_{l=i+1}^{r}
k_l,
&
i=j
\\
\\
\phi^{\dagger}_i
\phi_j
\prod_{l=i+1}^{r}
k_l,
&
i>j
\\
\\
0,
&
i<j
\end{array}
\right.
\quad\ \text{for }\ 
s \leq i \leq r,\ 
s+1 \leq j \leq r.
\end{align*}
%
This general rank solution for $x=1$ in terms of $t$-bosons was recently obtained in \cite{ProlhacEM, AritaAMP}. A generalisation of these results that includes a spectral parameter was found earlier in \cite{InoueKO} and independently in the case of super-algebras in \cite{Tsuboi}.

We also introduce a twist operator 
\begin{align*}
S^{(s)}=\prod_{l=s+1}^{r} k_l^{(l-s)u}
\ \ \text{for }\
1 \leq s \leq r-1,
\qquad
S^{(r)} = 1,
\end{align*} 
where we perform the reparametrization $q=t^{u}$. The operators  $\{k,\phi, \phi^{\dagger}\}$ are generators of the $t$-boson algebra \eqref{t-bos} with subscripts used to denote commuting copies of the algebra. Note that all operators in $\widetilde{L}^{(s)}(x)$ implicitly also carry an index $(s)$, as we will adopt the convention that operators in $\widetilde{L}^{(s)}(x)$ and $\widetilde{L}^{(s')}(x)$ (as well as $S^{(s)}$ and $S^{(s')}$) commute for $s \not= s'$.

\section{A polynomial example}
\label{se:polexample}
We look at an explicit example for rank 2 taking $\delta=(0,0,1,1,2,2)$. In this case $f_\delta$ corresponds to the nonsymmetric Macdonald polynomial $E_\delta$ \cite{CantinidGW}, which, using the notation $q=t^{u}$, $ [m]=(1-t^{m})/(1-t)$ is given by
\begin{multline}
E_\delta(x_1,\ldots,x_6;q,t)= x_3 x_4 x_5^2 x_6^2 + \frac{t^2}{[3+u]} (x_1+x_2)x_3x_4 x_5x_6(x_5+x_6) \\
+ \frac{t^{4}[2]}{[3+u][4+u]}x_1x_2x_3x_4x_5x_6.
\label{eq:Edeltar=2}
\end{multline}

We now verify the matrix product form 
\[
\Omega_{\delta^+} E_\delta(x_1,\ldots,x_6;q=t^{u},t)=\Tr \big[ A_0(x_1)A_0(x_2)A_1(x_3)A_1(x_4)A_2(x_5)A_2(x_6)S\big],
\]
for this explicit solution. From \eqref{eq:Arank2} we see that
\begin{align*}
A_0(x) = 1+x\phi,\quad A_1(x) &= x k,\quad A_2(x) = x \phi^\dag + x^2,
\end{align*}
and using \eqref{Exchange2} we note that $S$ should satisfy
\[
qS\phi - \phi S =0, \qquad S\phi^\dag -q \phi^\dag S =0.
\]

An explicit representation for the $t$-bosons in terms of semi-infinite matrices is given by 
\be
\label{eq:t-bosrep}
(\phi)_{m,m'} = (1-t^m) \delta_{m,m'-1},\quad (\phi^\dag)_{m,m'} = \delta_{m,m'+1},\quad (k)_{m,m'}=t^m \delta_{m,m'},
\ee
and $S$ has the form
\[
S = k^{u} = \diag\{1,t^{u},t^{2u},\ldots\} = \diag\{1,q,q^{2},\ldots\}.
\]
Up to a normalization, the nonsymmetric Macdonald polynomial $E_\delta$ is now represented in matrix product form by
\begin{align*}
&\Tr \left[ \left(1+ x_1 \phi\right) \left(1+ x_2 \phi\right) x_3 k x_4 k  x_5 \left(\phi^\dag+ x_5 \right) x_6 \left(\phi^\dag+ x_6 \right)S\right] \nonumber\\
&=x_3x_4x_5x_6 \Tr  \left[\left( x_5 x_6 k^2 + (x_1+x_2)(x_5+x_6)\phi k^2 \phi^\dag + x_1x_2 \phi^2 k^2 (\phi^\dag)^2 \right)S\right],
\end{align*}
where other terms involving unequal powers of $\phi$ and $\phi^\dag$ have zero trace.
Normalising with $\Omega_{\delta^+}=\Omega_{221100}=\Tr(k^2S)$ we finally get
\begin{multline*}
E_\delta(x_1,\ldots,x_6;q=t^{u},t) = x_3x_4x_5^2x_6^2 + x_3x_4x_5x_6(x_1+x_2)(x_5+x_6) t^2 
\frac{\Tr \phi \phi^\dag k^2S}{\Tr k^2S} \\ + x_1x_2x_3x_4 x_5x_6 t^4 \frac{\Tr \phi^2 (\phi^\dagger)^2 k^2S}{\Tr k^2S},
\end{multline*}
which can be shown to equal \eqref{eq:Edeltar=2} using the explicit representation \eqref{eq:t-bosrep}. 

\section{Conclusion}

We have derived a new explicit formula for symmetric Macdonald polynomials using the matrix product formalism. Our main results are the matrix product formula Theorem~\ref{th:MP} for the distinguished basis $\{f_{\sigma\circ\lambda}\}_{\sigma \in S_n}$ in the ring of homogeneous polynomials of degree $\lambda$. In the specialisation $q=1$ and $x_1=x_2=\dots = x_n=1$ these give the (unnormalized) stationary probabilities for the multi-species asymmetric exclusion process. The limits $t\rightarrow 0$ (or $t\rightarrow \infty$) give the stationary probablities for the totally asymmetric exclusion process, and a recent derivation of the matrix product formula for those was given in \cite{KunibaMO} using the tetrahedron equation for three-dimensional integrability. While there are several similarities with our approach, the resulting expressions for the matrices in terms of bosonic operators in \cite{KunibaMO} are different from those exhibited here.  

By Lemma~\ref{le:symmetry} the matrix product formula for the polynomials $f$ leads to the result for Macdonald polynomials. We mention here that the normalization factor $\Omega_\lambda$ in Theorem~\ref{th:MP} can be calculated explicitly,
\begin{align}
\label{norm}
\Omega_{\lambda}(q,t)
= 
\prod_{1 \leq i < j \leq r}
\frac{1}{1-q^{j-i} t^{\lambda'_i-\lambda'_j}},
\end{align}
where $r$ is the largest part of $\lambda$. As a nontrivial corollary we obtain a new summation formula for Macdonald polynomials, presented in Theorem~\ref{th:sum}. We give explicit examples of several of our constructions.

\acknowledgements
\label{sec:ack}
We thank the Galileo Galilei Institute and the organisers of the research program \textit{Statistical Mechanics, Integrability and Combinatorics} for kind hospitality during part of this work. It is a pleasure to thank Philippe Di Francesco for bringing the papers \cite{Kirillov_Noumi1,Kirillov_Noumi2} to our attention, and Eric Ragoucy, Ole Warnaar and Paul Zinn-Justin for helpful remarks and extended discussions on related topics. JdG and MW are generously supported by the Australian Research Council (ARC) and the ARC Centre of Excellence for Mathematical and Statistical Frontiers (ACEMS).


\end{document}